%Paper: alg-geom/9601023
%From: Tregrob@aol.com
%Date: Tue, 23 Jan 1996 08:17:34 -0500

\magnification =\magstep 2
\hoffset = -0.3 in
\hsize = 5 in
\voffset = -0.3 in
\vsize = 6.6 in

\mathsurround = 1pt
\baselineskip=13.5 pt

\centerline{\bf COMPACTIFICATIONS OF FAMILIES OF  CURVES}
\medskip
\centerline{Robert Treger}

\line{\empty}
\noindent
{\bf 1. Introduction and Statement of the Problem}
\medskip \noindent
The purpose
of this note is to describe a good compactification of
the variety $$V_{n,d}\subset {\bf P}^N \qquad (N=n\,(n+3)/2)$$
of irreducible
complex projective plane curves of degree $n$ with $d$
nodes and no other singularities. The problem is
discussed in classical literature, as well as
in %
[1, Sect.~1], [2, p.~219], [3, Sect.~2], etc. The naive
compactification of $V_{n,d}$, namely $\overline {V}_{n,d}\subset
 {\bf P}^N$, has very
complicated singularities.  In particular, $\overline {V}_{n,d}$ is not
unibranch. One can construct a unibranch
compactification as follows.  For $d \geq 1$, let
 $$\Sigma_{n,d}\subset
{\bf P}^N
\times {\rm Sym}^d({\bf P}^2)$$
be the closure of
the locus of pairs
$(E,\>
\sum_{i=1}^dP_i)$, where $E$ is an irreducible nodal curve
and $P_1 {\rm,\dots ,}P_d$ are its nodes.  We have proved that
$\Sigma_{n,d}$
 is unibranch everywhere [4]. However, $\Sigma_{n,d}$ is
not a normal variety, even though it has no singularities
in codimension $1$ [1]. We would like a
compactification to be a normal variety with at worst
rational singularities.  Our compactification is in fact
a quotient of a smooth variety by the full symmetric
group ${\cal S}_d$ (Theorem 2.8).

\line{\empty}
\noindent
 {\bf 2. Compactifications of Configuration Spaces and} ${\bf
V_{n,d}}$
\medskip \noindent
{\bf 2.1.} To begin with, we will describe a nice
compactification of ${\rm Sym}^d({\bf P}^2)\backslash \Delta$,
where $\Delta$ denotes the
singular locus. The construction presented below can be
generalized to any smooth projective variety of dimension
at least $2$ (in place of ${\bf P}^2$).

Let $d \geq 1$ be an
integer.
Let $L_s(P_1 {\rm,\dots ,}P_d)$
denote the linear system of
curves of degree $s$ having singularities in $d$ distinct
points $P_1 {\rm,\dots ,}P_d \in {\bf P}^2$, and
${\bf P}L_s(P_1 {\rm,\dots ,}P_d)
 \subset {\bf P}^{s(s+3) \! /2}$
the corresponding subspace.  Let $k = k(d)$  be
the smallest integer such that $\dim
L_k(\! P_1\! {\rm,\dots }, \! P_d \! )$ is independent of the choice of
$\! P_1\! {\rm,\dots }, \! P_d \! $.   We get a rank $k(k+3) \! /2 + 1
 - 3d$
vector subbundle
${\cal E}(d,k)$ of the trivial vector bundle on $
{\rm Sym}^d({\bf P}^2)\backslash \Delta$
whose fibers are linear systems of all plane curves of
degree $k$.
It can be extended to a reflexive sheaf
$\alpha _\ast {\cal E}(d,k)$ on ${\rm Sym}^d({\bf P}^2)$,
where $\alpha {\rm :}\,{\rm Sym}^d({\bf P}^2)\backslash \Delta
\rightarrow
{\rm Sym}^d({\bf P}^2)$ is the natural inclusion.
Clearly $\alpha _\ast {\cal E}(d,k)$
is not locally free provided $d \geq 2$.

 Let $G_{d,k} =
{\bf Gr}(k(k+3) \! /2+1-3d, \, k(k+3) \! /2+1)$ denote the
Grassmannian of $(k(k+3) \! /2 - 3d)\,$-\thinspace planes in ${\bf
P}^{k(k+3) \! /2}$.
We have a natural inclusion
$$\gamma_k \, {\rm :\ \ }  {\rm Sym}^d({\bf P}^2)\backslash \Delta
\rightarrow
G_{d,k}, \qquad \sum_{i=1}^dP_i \mapsto {\bf P}L_k(P_1 {\rm,\dots
,}P_d).$$

\bigskip \noindent
 {\bf 2.2. Definition.}
The variety $\Gamma_d = \Gamma_{d,k} = \overline{\gamma_k
({\rm Sym}^d({\bf P}^2)\backslash \Delta)}$ is said to
be a {\it compactification\/} of ${\rm Sym}^d({\bf P}^2)
\backslash \Delta$.

\bigskip \noindent
{\bf 2.3.}
Let
$F ({\bf P}^2,d)$ denote the complement in $({\bf P}^2)^d$ of the
large
diagonals. To a point $(P_1 {\rm,\dots ,} P_d)\in F({\bf P}^2,d)
$, we can
associate a flag

$${\bf P}L_k(P_1 {\rm,\dots ,}P_d) \subset \ldots \subset {\bf
P}L_k(P_1)
 \subset {\bf P}^{k(k+3) \! /2}.$$
 We obtain a natural inclusion
 $$\phi_k \, {\rm :\ \ } F({\bf P}^2,d) \rightarrow  F_{d,k}$$
where $F_{d,k} \subset G_{d,k} \,\times \ldots \times \,G_{1,k}$ is
an appropriate
flag manifold.

\bigskip \noindent
 {\bf 2.4. Definition.}
The variety $\Phi_d = \Phi_{d,k} = \overline{\phi_k
(F({\bf P}^2,d))}$ is said to
be a {\it compactification\/} of $F({\bf P}^2,d)$.

\bigskip \noindent
 {\bf 2.5.}
For each $r \geq k$, one can
similarly define $\Gamma_{d,r}$ and $\Phi_{d,r}$. We have two
natural morphisms
$$\Phi_{d,r}\rightarrow \Gamma_{d,r}\rightarrow {\rm Sym}^d({\bf
P}^2).$$

Consider the restriction of the tautological vector
subbundle over
$G_{d,k}$ to $\Gamma_d$. By the universal property
of Grassmannians, it is a unique minimal vector bundle
extension of ${\cal E} (d,k)$. Similarly, let
$${\cal T}_d \subset  \ldots \subset {\cal T}_1
\qquad ({\rm {rk}} \,{\cal T}_i = k(k+3) \! /2+1-3i).$$
denote the restriction to
$\Phi_d$ of the universal flag of vector bundles over $F_{d,k}.$

\bigskip \noindent
  {\bf 2.6.}
We fix $d \geq 1$ and vary $n$ and $g = (n-1)(n-2) \! /2-
d.$ Consider a natural inclusion
$$V_{n,d}\subset {\bf P}^N \times {\Gamma}_d$$
sending a
curve $C\in V_{n,d}$ to a pair $(C,\,{\bf P}L_k(P_1
{\rm,\dots ,}P_d))$,
where
$P_1 {\rm,\dots ,}P_d$ are the nodes of $C$.
Let $V^{'}_{n,d} \subset {\bf P}^N \times \Phi_d$
denote the inverse image of $V_{n,d}$ via the morphism
$\Phi_d \rightarrow \Gamma_d.$
Let
$${\bf E}(n,g)\subset {\bf P}^N \times \Gamma_d, \qquad
{\bf F}(n,g)\subset {\bf P}^N \times \Phi_d$$
be the closures of $V_{n,d}$ and $V^{'}_{n,d}$, respectively.
The group ${\cal S}_d$ acts on ${\bf F}(n,g)$ via its canonical action
on
$\Phi_d$. We have two natural morphisms
$${\bf F}(n,g)\rightarrow{\bf E}(n,g)\rightarrow \Sigma_{n,d}.$$

 For $n \geq k$, ${\bf F}(n,g)$ is just the projective
bundle ${\bf P}({\cal T}_d)$ over $\Phi_d$.
In general, ${\bf F}(n,g)$ is {\it not\/} a
projective bundle over the image of its projection to
$\Phi_d$.  We can describe boundary points of $\Phi_d$ and prove
the theorem below.

\medskip \noindent
 {\bf 2.7. Definition.} The ${\bf E} (n,g)$ is said to be
a {\it compactification\/} of $V_{n,d}$.
\medskip \noindent
\proclaim 2{.}8{.} Theorem.
{\rm i)} ${\bf F}(n,g)$ is
a smooth variety;
\hfil \break
\indent
 {\rm ii)} $ {\bf E}(n,g) = {\bf F}(n,g)/{\cal S}_d$ is a normal
variety with at worst rational singularities;
\hfil \break
\indent
 {\rm iii)} $\Phi_{d,r}$
is a smooth variety, and $\Phi_{d,r}/{\cal S}_d = \Gamma_{d,r}\>
(r \geq k(d))${\rm ;}
\hfil \break
\indent
{\rm iv)} $\Gamma_{d,k}$ is canonically isomorphic to
$\Gamma_{d,r}$,
and $\Gamma_{d,r}$
is a normal variety with at worst rational singularities
$(r \geq k(d))$.
\medskip \noindent
{\bf References}\par
\frenchspacing
\item{[1]}
S.~Diaz \& J.~Harris,
{\sl Geometry of the Severi varieties},
Trans.\ Amer.\ Math.\ Soc. {\bf 309} (1988) 1--34.
\item{[2]}
D.~Eisenbud \& J.~Harris,
{\sl Progress in the theory
of algebraic curves},
Bull.\ Amer.\ Math.\ Soc. {\bf 21} (1989) 205--232.
\item{[3]}
W.~Fulton,
{\sl On nodal curves},
Algebraic
geometry---Open problems, Proceedings, Ravello 1982,
C.~Ciliberto, F.~Chione and F.~Orecchi (eds.),
Lect.\ Notes
Math. {\bf 997}, Springer-Verlag, Berlin\ \thinspace New\thinspace
York, 1983,
pp.~146--155.
\item{[4]}
R.~Treger,
{\sl Local properties of families of plane curves}, J.~Differential
Geometry
{\bf 39} (1994) 51--55.
\medskip
\parindent = 2 in
AG{\&}N, Princeton, New Jersey

tregrob@aol.com
\vfill \eject
\bye